\documentclass[12pt]{article}
\usepackage{graphicx}
\usepackage{amsmath}
\usepackage{amsfonts}
\usepackage{amssymb}

\newtheorem{theorem}{Theorem}

\newtheorem{conjecture}[theorem]{Conjecture}
\newtheorem{corollary}[theorem]{Corollary}

\newtheorem{proposition}[theorem]{Proposition}
\newtheorem{remark}[theorem]{Remark}

\newenvironment{proof}[1][Proof]{\textbf{#1.} }{\ \rule{0.5em}{0.5em}}

\def\text{\hbox}

\def\S{{\bf S}}
\def\E{\mathcal{A}}
\def\C{\mathcal{C}}

\begin{document}

\title{Complexity, Heegaard diagrams and generalized Dunwoody manifolds \footnote{Work performed under the auspices
of the G.N.S.A.G.A. of I.N.d.A.M. (Italy) and the University of
Bologna, funds for selected research topics. The third author was
partially supported by the grant NSh-5682.2008.1, by the grant of
the RFBR, and by the grant of the Siberian Branch of RAN}}

\author{Alessia Cattabriga \and Michele Mulazzani \and Andrei Vesnin}

\maketitle

\begin{abstract}
We deal with Matveev complexity of compact orientable 3-manifolds represented via Heegaard diagrams.
This lead us to the definition of modified Heegaard complexity of  Heegaard diagrams and of  manifolds. We define a class of manifolds
which are generalizations of Dunwoody
manifolds, including cyclic branched coverings of two-bridge knots
and links, torus knots, some pretzel knots, and some theta-graphs.
Using modified Heegaard complexity,  we
obtain upper bounds for their Matveev complexity, which linearly
depend on the order of the covering. Moreover,  using homology arguments due to Matveev and Pervova we obtain lower bounds.
\end{abstract}

\maketitle

\section{Introduction and preliminaries}

The notion of complexity for compact 3-dimensional manifolds has been introduced by S. Matveev via simple spines.
We briefly recall its definition (for further reference see
\cite{M1,M2}).

A polyhedron $P$ embedded into a compact connected 3-manifold $M$ is
called a \textit{spine} of  $M$ if $M$ collapses to $P$ in the case
$\partial M \neq \emptyset$, and if $M-{\rm Int}(B)$ collapses to
$P$ in the case $\partial M = \emptyset$, where $B$ is a closed 3-ball in
$M$. Moreover, a spine $S$ is said to be \textit{almost simple} if
the link of each point $x\in S$ can be embedded into $K_4$, a
complete graph with four vertices. A \textit{true vertex} of an
almost simple spine $S$ is  a point $x\in S$ whose link is
homeomorphic to $K_4$.

The \textit{complexity} $c(M)$ of $M$ is the minimum number of true
vertices among all almost simple spines of $M$. Complexity is
additive under connected sum of manifolds and, for any
integer  $n \geqslant 0$, there are only finitely many closed prime
manifolds with complexity $n$.

In the closed orientable case there are  only  four prime manifolds
of complexity zero which are  $\S^3$, $\mathbb{RP}^3$, $\S^2\times
\S^1$, and $L_{3,1}$. Apart from these special cases, it can be
proved that $c(M)$ is the minimum number of tetrahedra needed to
obtain $M$ by pasting together their faces (via  face paring).  A
complete classification of closed orientable prime manifolds up to
complexity 12 can be found in \cite{M3,M4}.

In general, the computation of the complexity of a given manifold is a
difficult problem. So, two-sided estimates of complexity become
important, especially when dealing with infinite families of
manifolds (see, for example, \cite{M2,MPV,PV}).

By \cite[Theorem 2.6.2]{M2},  a lower bound  for the complexity of a
given manifold can be obtained via the computation of its first
homology group.  Moreover, for a hyperbolic
manifold a lower bound can be obtained via volume arguments  (see
\cite{M2,MPV,PV}). On the other hand, upper bound can be found using
triangulations.

In this paper we deal with the possibility of calculating complexity
via Heegaard decompositions. This way of representing 3-manifold has
revealed to be very useful in different contests. So, it is natural
to wonder whereas it is possible to calculate complexity via
Heegaard diagrams. In Section~\ref{sec2} we use Heegaard diagrams to
define modified Heegaard complexity of  compact 3-manifolds and
compare this notion with Matveev complexity. A widely studied
family of manifolds, defined via Heegaard diagrams, is the one of
Dunwoody manifolds (see \cite{Du}). This family  coincides with the
class of strongly-cyclic branched coverings of $(1,1)$-knots (see
\cite{CM1}), including, for example, 2-bridge knots, torus knots and
some  pretzel knots. In  Section~\ref{sec3} we construct a  class of
manifolds that  generalizes the class of Dunwoody manifolds,
including  other   interesting class of manifolds such as
cyclic-branched coverings of 2-component 2-bridge links. In
Section~\ref{sec4}, using modified Heegaard complexity, we obtain
 two-sided estimates for the
complexity of some families of generalized Dunwoody manifolds.

\section{Heegaard diagrams and complexity} \label{sec2}

In this section we introduce the notions of modified complexity for
Heegaard diagrams and for  manifolds, comparing these notions with  Matveev
complexity of manifolds. Let us start by recalling some
definitions.

Let $\Sigma_g$ be a closed, connected, orientable surface of genus
$g$.  A \textit{system of curves} on $\Sigma_g$ is a (possibly
empty) set of simple closed curves $\C=\{\gamma_1,\ldots,\gamma_k\}$
on $\Sigma_g$ such that $\gamma_i \cap \gamma_j = \emptyset$ if
$i\ne j$, for $i,j=1,\ldots, k$. Moreover, we denote with
$V(\mathcal{C})$  the set of connected components of the surface
obtained by cutting $\Sigma_g$  along the curves of $\mathcal{C}$.
The system $\mathcal{C}$ is said to be \textit{proper} if all
elements of $V(\mathcal{C})$ have genus zero, and \textit{reduced}
if  either $\vert
V(\mathcal{C})\vert =1$ or $V(\C)$ has no elements of genus zero. Thus, $\mathcal{C}$ is: (i) proper and
reduced if and only if it consists of one element of genus $0$; (ii)
non-proper and reduced if and only if all its elements are of genus
$>0$; (iii) proper and non-reduced if and only if it has more than one element and all of them are of genus
$0$; (iv) non-proper and non-reduced if and only if it has at least
one element of genus $0$ and at least one element of genus $>0$.
Note that a proper reduced system of curves on $\Sigma_g$ contains
exactly $g$ curves.

We denote by $G(\mathcal C)$ the graph which is
dual to the one determined by $\mathcal{C}$ on $\Sigma_g$. Thus,
vertices of $G(\mathcal C)$ correspond to elements of
$V(\mathcal{C})$ and edges correspond to curves of $\mathcal{C}$.
Note that loops and multiple edges may arise in~$G(\mathcal C)$.

A \textit{compression body} $K_g$ of genus $g$ is a 3-manifold with
boundary, obtained from $\Sigma_g\times [0,1]$ by attaching a finite
set of 2-handles $Y_1,\ldots, Y_k$ along a system of  curves (called
\textit{attaching circles}) on $\Sigma_g\times\{0\}$ and filling in
with balls all the spherical boundary components of  the resulting
manifold, except from $\Sigma_g\times\{1\}$ when $g=0$. Moreover,
$\partial_+ K_g=\Sigma_g\times\{1\}$ is called the \textit{positive}
boundary of $K_g$, while $\partial_- K_g = \partial K_g-\partial_+
K_g$ is called \textit{negative} boundary of $K_g$. Notice that a
compression body is a handlebody if an only if $\partial_- K_g =
\emptyset$, i.e., the system of the attaching circles on
$\Sigma_g\times\{0\}$ is proper. Obviously homeomorphic compression
bodies  can be obtained with
(infinitely many) non isotopic systems of attaching circles.

\begin{remark} \label{riduzione}
\textup{If the system of attaching circles is not reduced then it
contains at least one reduced subsystem of curves  determining the
same compression body $K_g$. Indeed, if $\C$ is the system of
attaching circles, denote with  $V^+(\mathcal{C})$  the set of
vertices of $G(\C)$ corresponding to the components with genus
greater then zero, and  with $\mathcal{A}(\mathcal{C})$ the set
consisting of  all the graphs $T_i$ such that: 
\begin{itemize}
\item $T_i$ is a subgraph of $G(\mathcal C)$;
\item if $V^+(\mathcal{C})=\emptyset$ then $T_i$ is a maximal tree in
$G(\mathcal C)$;
\item if $V^+(\mathcal{C})\ne \emptyset$ then $T_i$ contains all the
vertex of $G(\mathcal C)$ and each component of $T_i$ is a tree
containing exactly a vertex of $V^+(\mathcal{C})$.
\end{itemize}
Then, for any  $T_i \in \mathcal{A}(\mathcal{C})$, the system of curves
obtained by removing from $\C$ the curves corresponding to the edges
of $T_i$ is reduced and determines the same compression body. Note
that this operation corresponds to removing complementary 2- and
3-handles. Moreover, it is easy to see that  if $\partial_- K_g$ has
$k$ boundary components with genus $g_1,\ldots,g_k$ then
$$
\vert E(T_i)\vert = \vert\C\vert - n - k + 1 + \sum_{j=1}^k
g_j
$$
for each $T_i\in\mathcal{A}(\C)$, where $E(T_i)$ denotes the edge set of
$T_i$.}
\end{remark}

Let $M$ be  a compact, connected, orientable 3-manifold without
spherical boundary components. A \textit{Heegaard surface} of genus
$g$ for $M$ is a surface $\Sigma_g$  embedded in $M$ such that
$M-\Sigma_g$ consists of two components whose closures $K'$ and
$K''$ are (homeomorphic to), respectively, a genus $g$ handlebody
and  a genus $g$ compression body.

The triple $(\Sigma_g, K',K'')$ is called \textit{Heegaard
splitting} of $M$. It is a well known
fact that each compact connected orientable 3-manifold without
spherical boundary components admits a Heegaard splitting.

\begin{remark} \textup{By \cite[Proposition 2.1.5]{M2} the complexity of a manifold is not affected by puncturing it. So, with the aim of computing complexity, it is not restrictive assuming that the manifold does not have spherical boundary components.}
\end{remark}

On the other hand, a triple $H=(\Sigma_g,
\mathcal{C'},\mathcal{C''})$, where  $\mathcal{C'}$ and
$\mathcal{C''}$ are two systems of curves on $\Sigma_g$, such that
they intersect transversally and $\mathcal{C'}$ is proper, uniquely
determines a 3-manifold $M_H$ which corresponds to the Heegaard
splitting $(\Sigma_g, K',K'')$, where $K'$ and $K''$ are respectively
the handlebody and the compression body whose attaching circles
correspond to the curves in the two systems. Such a triple is called
\textit{Heegaard diagram}  for $M_H$.

We denote by $\Gamma(H)$  the graph embedded in $\Sigma_g$, obtained
from the curves of $\mathcal{C'}\cup \mathcal{C''}$, and by
$\mathcal{R}(H)$ the set of regions of $\Gamma(H)$. Note that
$\Gamma(H)$ has two types of vertices: singular vertices which are
4-valent and non-singular ones which are 2-valent. A diagram $H$ is
called \textit{reduced} if both the systems of curves are reduced.
If $H$ is non-reduced, then we denote by $\textup{Rd}(H)$ the set of
reduced Heegaard diagrams obtained from $H$ by reducing the system
of curves.

In \cite[Section~7.6]{M2} the notion of complexity of a reduced
Heegaard diagram $H$ of a genus two closed manifold is defined as
the number $c(H)$ of singular vertices of the graph $\Gamma(H)$.
Moreover the author proved that $c(M_H) \leqslant c(H)$.

Now we extend this definition to the general case,  slightly
modifying it  in order to obtain a better estimate for the
complexity of $M_H$.

The modified complexity of  a reduced Heegaard diagram  $H$ is
$$\widetilde{c}(H)=c(H)-\max\,\{n(R)\mid R\in\mathcal{R}(H)\},$$ where
$n(R)$ denotes the number of singular vertices contained in the
region $R$, and the modified complexity of a (non-reduced) Heegaard
diagram $H$ is $$\widetilde{c}(H)=\min\,\{\tilde{c}(H')\mid
H'\in\textup{Rd}(H)\}.$$

We define the \textit{modified Heegaard complexity} of a closed
connected 3-manifold $M$ as
$$
\widetilde{c}(M) = \min\,\{\widetilde{c}(H) \mid H
\in\mathcal{H}(M)\} ,
$$
where $\mathcal{H}(M)$ is the set of all Heegaard diagrams of $M$.

The following statement generalizes a result of \cite[Proposition
2.1.8]{M2} (for the case of reduced diagrams of closed manifolds)
and \cite{C} (for case of Heegaard diagrams arising from gem
representation of closed manifolds).

\begin{proposition} \label{prop4}
If $M$ is a compact connected 3-manifold then
$$
c(M) \leqslant \widetilde{c}(M) .
$$
\end{proposition}

\begin{proof}
Let $H=(\Sigma_g,\mathcal C',\mathcal C'')$ be a Heegaard diagram of
$M$ and let $(\Sigma_g,K',K'')$ be the associated Heegaard
splitting. We want to prove that $c(M) \leqslant \widetilde c(H)$.
From the definition of modified complexity it is clear that we can
suppose that $H$ is reduced. If $\partial M = \emptyset$ then the
statement is given in \cite[Proposition 2.1.8]{M2}. For the case
$\partial M\ne \emptyset$ the same proof works because of the
following reason. The simple polyhedron obtained as the union of
$\Sigma_g$ with the core of the 2-handles of $K'$ and $K''$ is a
spine with $c(H)$ singular vertices of $M-{\rm Int}(B)$, where $B \subset
K'$ is a closed ball. Since $\partial M$ is contained in $K''$, a
spine for $M$ can be obtained by connecting $\partial B$ with
$\partial M$ via pinching a region of~$\mathcal{R}(H)$.
\end{proof}

%\bigskip

By results of \cite{CC},  the upper  bound in Proposition~\ref{prop4}
becomes an equality for the $69$ closed connected prime orientable
3-manifolds admitting a (colored) triangulation with at most $28$
tetrahedra. As far as we know there is no example where the
strict inequality holds.

\begin{conjecture}
For every compact connected orientable 3-manifold $M$ the equality
$c(M)= \widetilde{c}(M)$ holds.
\end{conjecture}

\section{Generalized Dunwoody manifolds} \label{sec3}
In this section we define a class of manifolds that generalizes the
class of Dunwoody manifolds introduced in \cite{Du}.

\begin{figure}
 \begin{center}
 \includegraphics*[totalheight=6.5cm]{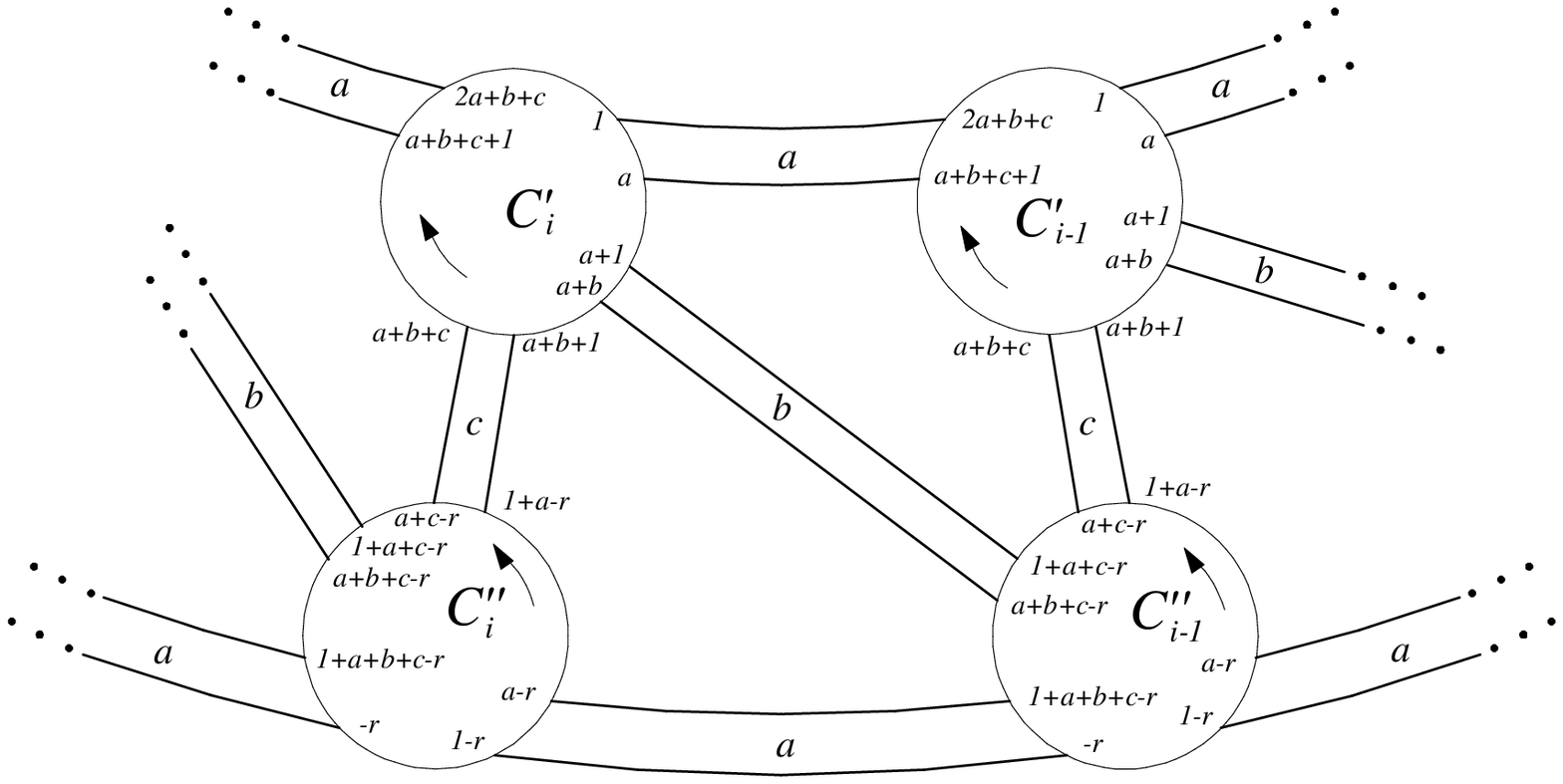}
 \end{center}
 \caption{A Dunwoody diagram.}
 \label{dun}
\end{figure}

A \textit{Dunwoody diagram} is a trivalent regular planar graph,
depending on six integers $a,b,c,n,r,s$, such that $n>0$, $a, b ,c
\geqslant 0$ and $d=2a+b+c>0$, and it is defined as follows (see
Figure~\ref{dun}).

It contains $n$ internal circles $C'_1,\ldots,C'_n$, and $n$
external circles $C''_1,\ldots,C''_n$, each having $d$ vertices. The
circle $C'_i$ (resp. $C''_i$) is connected to the circle $C'_{i+1}$
(resp. $C''_{i+1}$) by $a$ parallel arcs, to the circle $C''_{i}$ by
$c$ parallel arcs and to the circle $C''_{i-1}$ by $b$ parallel
arcs, for every $i=1,\ldots,n$ (subscripts mod $n$). We denote by
$\E$ the set of arcs, and by $\mathcal B$  the set of circles. By
gluing the circle $C'_i$ to the circle $C''_{i+s}$ in the way that
equally labelled vertices are identified together  (see
Figure~\ref{dun} for the labelling), we obtain a Heegaard diagram
$H(a,b,c,n,r,s)=(\Sigma_n,\C',\C'')$, where $\C'$ is the proper,
reduced system of curves arising from $\mathcal B$,
containing $n$ curves, and $\C''$ is the system of curves arising
from $\E$, containing $m>0$ curves.  Observe that the parameters $r$
and $s$ can be considered mod~$d$ and mod~$n$ respectively. We call
$H(a,b,c,r,n,s)$ \textit{closed Dunwoody diagram}. The
\textit{generalized Dunwoody manifold} $M(a,b,c,n,r,s)$ is the
manifold $M_{H(a,b,c,n,r,s)}$.

Since both the diagram and the identification rule are invariant
with respect to an obvious cyclic action of order $n$, the
generalized Dunwoody manifold $M(a,b,c,n,r,s)$ admits a cyclic
symmetry of order~$n$.

\begin{remark}\label{equivalenza}
\textup{It is easy to observe that diagrams $H(a,b,c,r,n,s)$ and
$H(a,c,b,d-r,n,n-s-1)$ are isomorphic, so they represent the same
manifold.  }
\end{remark}

A generalized Dunwoody manifold  $M(a,b,c,n,r,s)$ is a Dunwoody
manifold when the system $\C''$ of curves arising from $\E$ is
proper and reduced.  In this case $H(a,b,c,n,r,s)$ is a ``classical''
Heegaard diagram (see \cite{He}) and therefore all Dunwoody
manifolds are closed.

As proved in \cite{CM}, the class of Dunwoody manifolds coincides
with the class of strongly-cyclic branched  covering of
$(1,1)$-knots. So, in particular, it contains all cyclic branched
coverings of 2-bridge knots. It is not known if cyclic branched
coverings of 2-bridge links (with two components) admit
representations as Dunwoody manifolds, but they surely are
generalized Dunwoody manifolds. This can be shown by introducing a
polyhedral description for generalized Dunwoody manifolds.

\begin{figure}
\begin{center}
\includegraphics*[totalheight=10cm]{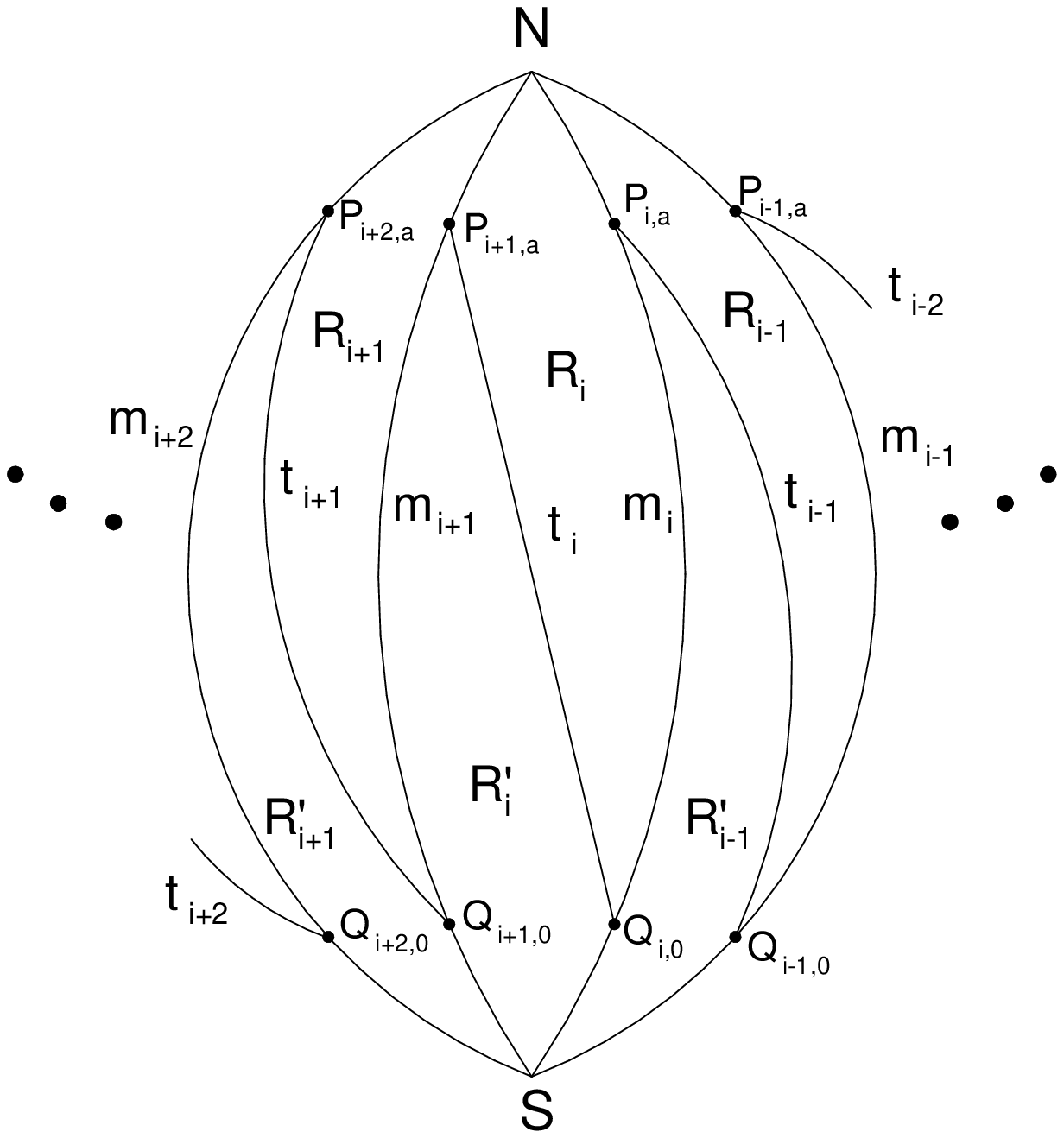}
\end{center}
\caption{Polyhedral description of generalized Dunwoody manifolds.}
\label{pol}
\end{figure}

Referring to Figure \ref{pol}, let $B$ be the closed unitary  3-ball in
$\mathbb{R}^3$ and consider on its boundary $n$ equally spaced
meridians $m_1, \ldots m_n$ joining the north pole $N=(0,0,1)$ with
the south pole $S=(0,0,-1)$.  Subdivide each meridian $m_i$ into
$2a+b$ arcs with endpoints $P_{i,j}$, $j=0,\ldots,2a+b$, such that
$P_{i,0}=N$ and $P_{i,2a+b}=S$. Let $t_i \in \partial B$ be the
shortest arc connecting  $P_{i,a+b}$ with $P_{i+1,a}$, for
$i=1,\ldots,n$. We subdivide $t_i$ into $c$ arcs  with endpoints
$Q_{i,j}$ for $j=0,\ldots,c$ such that $Q_{i,0}=P_{i,a+b}$ and
$Q_{i,c}=P_{i+1,a}$. In this way $\partial B$ is subdivided into
$2n$ $d$-gons with $d=2a+b+c$. We denote by $R_1,\ldots,R_n$ the
$d$-gons containing  the north pole $P_{i,0} = N$ and by
$R'_1,\ldots,R'_n$ the $d$-gons containing the south pole. Moreover,
let
$$
P'_{i,0} = \begin{cases} P_{i,2a+b-r} & 0 \leqslant r \leqslant a, \\
Q_{i,r-a} & a \leqslant r \leqslant a+c, \\
P_{i+1,r-c} & a+c \leqslant r \leqslant 2a+b+c .
\end{cases}
$$
According to this definition $P'_{i,0}$ is a point on the boundary
of $R'_i$ obtained from $S$ by giving a combinatorial $r$-twist in
counterclockwise direction to the region $R'_i$.

We glue $R_i$ with $R'_{i+s}$ by an orientation reversing
homeomorphism matching the vertices of $R_i$ with the ones of
$R'_{i+s}$ such that $P_{i,0} \in R_i$ is identified with
$P'_{i+s,0} \in R'_{i+s}$. In this way we obtain a closed connected
orientable pseudomanifold $\widehat M(a,b,c,n,r,s)$ with a finite
number of singular points whose stars are cones over closed
connected orientable surfaces. By removing the interior of a regular
neighboorhood of each singular point we get a compact connected
orientable 3-manifold with (possibly empty) non-spherical boundary
components, which is homeomorphic to the generalized Dunwoody
manifold $M(a,b,c,n,r,s)$.

As a particular case,  an $n$-fold cyclic branched covering of a
2-bridge link/knot $\mathbf{b}(\alpha,\beta)$ is
$M(\beta,\alpha-2\beta,1,n,2\beta+1,s)$ where $s=(-1)^{\beta}$ if
$\mathbf{b}(\alpha,\beta)$ is a knot (i.e. $\alpha$ is odd) and
$s\ne 0$ if  $\mathbf{b}(\alpha,\beta)$ has two components (i. e.
$\alpha$ is even) (see \cite{Mi,Mu}).

\section{Upper and lower bounds} \label{sec4}
In this section we calculate the modified complexity of a closed
Dunwoody diagram in order to find upper bounds for the complexity of
some families of generalized Dunwoody manifolds.  For $n=1$, the
generalized  Dunwoody manifold is a a lens space (including
$\S^2\times \S^1$ and $\S^3$) in the closed case and a solid torus in
the case with boundary. Since the complexity of these manifolds has been
already studied (see \cite[Section 2.3.3]{M2}), we will always
suppose $n>1$.

\begin{theorem} \label{prop}
Let $H=H(a,b,c,n,r,s)=(\Sigma_n,\mathcal{C',\mathcal C''})$ be a
closed Dunwoody diagram, and $d=2a+b+c$. For each $\gamma \in \C''$
define $n(\gamma)$ as the number of singular vertices contained in
the cycle determined by $\gamma$ in $\Gamma(H)$. Then, with the
notation of Remark~\ref{riduzione} we have:
\begin{equation*} \label{formula}
\widetilde{c}(H) = nd - \max\left\{n(R)+\sum_{\gamma \in E(T)}
n(\gamma)\mid T\in\mathcal{A}(\C''), R\in\mathcal{R}(H_T)\right\},
\end{equation*}
where $E(T)$ is the edge set of the graph $T$ and $H_T$ is the
element of $\textup{Rd}(H)$ obtained by removing from $\C''$ the
curves belonging to $T$.
\end{theorem}

\begin{proof}
By construction the system $\C'$ is proper and reduced. The
statement follows from the definition of modified complexity and
Remark~\ref{riduzione}.
\end{proof}

This result allows us to find upper bounds for the modified complexity (and so for Matveev complexity) of generalized Dunwoody  manifolds. In the following subsections  we specialize the estimates  to the cases of
some important families.

\subsection{Dunwoody manifolds}

\begin{proposition} \label{prop8}
Let $M=M(a,b,c,n,r,s)$ be a Dunwoody manifold. Then
\begin{itemize}
\item[\textup{(i)}] If  $abc>0$ then
$$
c(M) \leqslant \left\{\begin{array}{ll} n(2a+b+c) -\max(2n,6) &
\textup{ if }
r\ne -b,-b\pm 1,\\
n(2a+b+c)-\max(2n,5) & \textup{ if }  r=-b\pm 1.\\
\end{array}\right.
$$
\item[\textup{(ii)}] If $abc=0$ and $\min(a,b+c)=0$ then
$$
c(M) \leqslant \left\{\begin{array}{ll} n(2a+b+c-4) & \textup{ if }
r \ne -b,-b\pm 1,\\
n(2a+b+c-3) & \textup{ if }  r=-b\pm 1.\\
\end{array}\right.
$$
\item[\textup{(iii)}] If  $abc=0$ and $\min(a,b+c)>0$ then
$$
c(M) \leqslant \left\{\begin{array}{ll} n(2a+b+c-2) & \textup{ if }
n > 3,\\
n(2a+c) -\max(2n,8-2k_0) & \textup{ if } n=2,3 , \,
b=0 \textup{ and } s=0,\\
n(2a+b) -\max(2n,8-k_0-k_1) & \textup{ if }  n=2, c=0 \textup{ and } s=0,\\
n(2a+b) -\max(2n,8-k_0) & \textup{ if }  n=3, c=0 \textup{ and } s=0, \\
n(2a+b) -\max(2n,8-k_1) & \textup{ if }  n=3, c=0 \textup{ and } s=1,\\
\end{array}\right.
$$
where $k_i=\left\{\begin{array}{ll} 2 & \textup{ if }  r=(-1)^i b,\\
1 & \textup{ if }  r=(-1)^i b\pm 1,\\ 0 & \textup{ otherwise}.
\end{array}\right.$
\end{itemize}
The cases not covered by the above formulas follow from the
homeomorphisms $M(a,b,c,r,n,s)\cong M(a,c,b,d-r,n,n-s-1)$ (see
Remark~\ref{equivalenza}).
\end{proposition}

\begin{proof} The graph $\Gamma(H)$ associated to a Heegaard diagram $H$
for $M(a,b,c,n,r,s)$ is obtained from the diagram depicted in
Figure~\ref{dun} by performing the prescribed identifications. Since
$\Gamma(H)$ is proper and reduced, then $G(\mathcal C'')$ is an
$n$-circle bouquet, so $T$ is a single point and therefore
$E(T)=\emptyset$. Hence by Theorem~\ref{prop}
$$
\widetilde{c}(H)\leqslant n(2a+b+c)-\max\{n(R)\mid
R\in\mathcal{R}(H)\}.
$$
In case (i) the upper (and lower) region of the Dunwoody diagram has
$2n$ vertices that are not identified together by the gluing, while
for all the other regions it is clear that $n(R)\leqslant 6$. More
precisely, the six vertices of hexagonal regions remain all distinct
if $r\ne-b,-b\pm 1$, while two of them are identified if $r=-b\pm
1$. If $r=-b$ then $M(a,b,c,n,r,s)$ is not a Dunwoody manifold since
$\Gamma(H)$  is not reduced.

In case (ii) the Dunwoody diagram has regions with $4n$ vertices. As
before, they remain all distinct under identifications if
$r\ne-b,-b\pm 1$, they become $3n$ if $r=-b\pm 1$, while if $r=-b$
the associated manifold is not Dunwoody.

In case (iii), if $n\geqslant 4$ then the upper (or lower) region
has $2n$ vertices while all other regions have at most $8$ vertices.
When $n=2$ or $n=3$ the computation is more tricky. We always have a
region with eight vertices, but, as before, some of them can be
identified together. Given such a maximal region, the number $k_i$
counts how many vertices of the circle $C'_i$ are identified with
the   ones of the circle $C''_{i+s}$.
\end{proof}

Proposition~\ref{prop8} allows to obtain an upper bound for the
complexity  of cyclic branched coverings of 2-bridge knots
(Corollary~\ref{Cor9}) and some families of torus knots
(Corollary~\ref{Cor11}), and a family of Seifert manifolds
(Corollary~\ref{Cor12}).

We recall that  $\mathbf{b}(\alpha,\beta)$ is a 2-bridge knot if and
only if $\alpha$ is odd.

\begin{corollary} \label{Cor9}
Let $C_n(\alpha,\beta)$ be the $n$-fold cyclic branched covering of
the 2-bridge knot  $\mathbf b(\alpha,\beta)$. Then for $n>2$ we have
$$c(C_n(\alpha,\beta)) \leqslant n(\alpha-2).$$
\end{corollary}

\begin{proof}
Since $\mathbf b(\alpha,\alpha-\beta)$ is the mirror image of
$\mathbf b(\alpha,\beta)$ we can suppose that $\beta$ is even. By
\cite{GM} we have that
$C_n(\alpha,\beta)=M((\alpha-1)/2,0,1,n,\beta/2,s)$, for a certain
$s=s(\alpha,\beta)$.
\end{proof}

This result improves the upper bound obtained in \cite{PV}, where
the lower bound has been obtained in the hyperbolic case (i.e.
$\beta\ne 1, \alpha-1$) via volume estimates. Now we give a lower
bound for the remaining cases.

\begin{proposition}
Let $n>2$. We have
$$
c(C_n(\alpha,1))=c(C_n(\alpha,\alpha-1)) \geqslant
\left\{\begin{array}{ll}
2 \log_5(\alpha/d)+d-2 & \textup{ if } n \textup{ is even,}\\
2(d-1)\log_5 2-1 & \textup{ if } n \textup{ is
odd,}\end{array}\right.
$$
where $d=\gcd (\alpha,n)$.
\end{proposition}

\begin{proof}
Obviously $C_n(\alpha,\alpha-1) \cong C_n(\alpha,1)$ since $\mathbf
b(\alpha,\alpha-1)$ is the mirror image of $\mathbf b(\alpha,1)$.
Moreover $\mathbf b(\alpha,1)$ is the torus knot of type
$(\alpha,2)$ and therefore $C_n(\alpha,1)$ is the Brieskorn manifold
of type $(2,\alpha,n)$ \cite{Mil}. Its first homology group is
$\mathbb{Z}^{d-1}\oplus \mathbb{Z}_{n/d}$ if $n$ is even, and
$\mathbb{Z}_2^{d-1}$ if $n$ is odd (see \cite{Ra, Ca}). Since the
manifold is irreducible (and different from $L_{3,1}$), the result follows by applying
Theorem~2.6.2 of \cite{M2}.
\end{proof}

\begin{corollary} \label{Cor11}
Let $T_n(k,h)$ be the $n$-fold cyclic branched covering of the torus
knot of type $(k,h)$. Then  we have
\begin{enumerate}
\item $c(T_n(k,h)) \leqslant n \, (2qk-2q-1)$ if $h=qk+1$ for $q>0$ and
$k>1$;
\item $c(T_n(k,h)) \leqslant n \, (2qk-2q-3)$ if $h=qk-1$ for $q, k > 1$;
\item $c(T_n(k,h)) \leqslant n \, (2q_1(s-1)(qq_1+1)+2qq_1-1)$
if $k=sq_1+1$ and $h=qk+s$ for $q, q_1 > 0$ and $s>1$.
\end{enumerate}
\end{corollary}

\begin{proof}
By \cite{AGM} we have  that
$$
T_n(k,qk+1)=M(1,k-2,(k-1)(2q-1),n,k,k)
$$
and
$$
T_n(k,qk-1)=M(1,k-2,(k-1)(2q-1)-2,n,(k-1)(2q-3),k).
$$
Moreover, by \cite{CM1}, there exists $s\in\mathbb{Z}$ such that
$$
T_n(sq_1+1, (sq_1+1)q+s) = $$
$$
= M(q_1,q_1(2qq_1(s-1)+2q+s-2),1+(s-2)q_1,2q_1^2(s-1)+sq_1+1).
$$
The result follows from Proposition~\ref{prop8}.
\end{proof}

We remark that an algorithm developed in \cite{CM1} allows us to
obtain a presentation of each $n$-fold cyclic branched covering of a
torus knot as a Dunwoody manifold and so to compute an upper bound
for the complexity by using Proposition~\ref{prop8}.

It is proved in \cite{GM2} that if $p > q > 0$ and $\gcd(p,q)=1$,
$n>1$,  $\ell > 0$, then Seifert manifolds
$$
S_n(p,q,\ell)= \{Oo,0\mid -1; \underbrace{(p,q), \ldots,
(p,q)}_{n-\textup{times}},(\ell, \ell-1)\}
$$
are Dunwoody manifolds that generalize the class of Neuwirth manifolds
introduced in \cite{Ne} and corresponding to $p=2$ and $q=\ell=1$.
Below we will give upper and lower estimates for complexity of these
Seifert manifolds.

\begin{corollary} \label{Cor12}
Suppose $\ell>1$ when $n=2$. The following estimate holds:
$$
c(S_n(p,q,\ell)) \leqslant n(p + q (n\ell - 2) -2).
$$
\end{corollary}

\begin{proof}
By results of \cite{GM2}, we have that
$$
S_n(p,q,\ell) = M(q,q(n\ell-2),p-2q,n,p-q,0)
$$
if $p \geqslant 2q$ and
$$
S_n(p,q,\ell) = M(p-q,2q-p,q(n\ell-2),n,p-q,1)
$$
otherwise. The result follows from Proposition~\ref{prop8}.
\end{proof}

\begin{proposition}
The following estimate holds:
$$
c(S_n(p,q,\ell)) \geqslant 2(n-1) \log_5p + 2
\log_5((n-1) \ell q-p)-1.
$$
\end{proposition}

\begin{proof}
Following \cite{Or}, a standard presentation of
$\pi_1(S_n(p,q,\ell))$ is
$$
\langle y_1,\ldots,y_n,y,h \mid [y_i,h], \, [y,h], \, y_i^p h^q, \,
y^{\ell} h^{\ell-1}, \, y_1 \cdots y_n y h; i=1,\ldots,n \rangle.
$$
By abelianization, we find that a presentation matrix for
$H_1(S_n(p,q,\ell))$ as a $\mathbb{Z}$-module is the circulant
matrix whose first row is given by the coefficient of $f(t) = - p +
\ell q \sum_{i=1}^{n-1}t^i$. By the theory of circulant matrices
\cite{BaM}, there exists a complex unitary matrix $F$, called
\textit{Fourier matrix}, such that
$$
FBF^* = D = \textup{Diag}(f(\zeta_1), f(\zeta_2), \ldots,
f(\zeta_n)),
$$
where $\zeta_1,\zeta_2,\ldots,\zeta_n$ are the $n$-roots of the
unity. So it follows that
$$
\vert Tor(H_1(S_n(p,q,\ell)))\vert = p^{n-1}((n-1) \ell q - p).
$$
Moreover, since $S_n(p,q,\ell)$ is irreducible (and different from
$L_{3,1}$), the result follows from Theorem~2.6.2 of \cite{M2}.
\end{proof}

\subsection{Cyclic branched coverings of two-bridge links}

We recall that  $\mathbf{b}(\alpha,\beta)$ is a 2-component 2-bridge
link if and only if $\alpha$ is even. In the next statement we deal with cyclic
branched coverings of  2-component 2-bridge
links of singly type (see \cite{MM}).
\begin{proposition}
Let $\mathbf{b}(\alpha,\beta)$ be a 2-bridge link with two
components and denote by $m_1$ and $m_2$ the homology classes of the
meridian loops of the two components. If $C_{n,s}(\alpha,\beta)$ is
the $n$-fold cyclic branched covering of $\mathbf{b}(\alpha,\beta)$
with monodromy $\omega(m_1)=1$, $\omega(m_2)=s \in \mathbb{Z}_n - \{
0 \}$ then
$$
c(C_{n,s}(\alpha,\beta)) \leqslant n(\alpha-2)+\frac{n}{d}-\alpha,
$$
where $d=\gcd(n,s)$.
\end{proposition}

\begin{proof}
By results of \cite{Mi,Mu}, we have
$$
C_{n,s}(\alpha,\beta)=M(\beta,\alpha-2\beta,1,n,2\beta+1,s),
$$
so we can use Theorem~\ref{prop} to calculate $\widetilde{c}(H)$
in order to obtain an upper bound for $c(C_{n,s}(\alpha,\beta))$.
The system of curves $\mathcal C''$ of the Dunwoody diagram
$$
H=H(\beta,\alpha-2\beta,1,n,2\beta+1,s)=(\Sigma_n,\mathcal{C'},\mathcal{C''})
$$
is not reduced. Indeed, taking advantage of its symmetries, it is
easy to see that it  consists of $n+d$ curves. More precisely, $d$
curves (that we call of type A) arise from all $n$ ``radial'' arcs
(i.e. the ones connecting the circles $C'_i$ and $C''_i$).
Each of these curves intersects $\mathcal{C'}$ in $n/d$ points. The
other $n$ curves (that we call of type B) arise from the remaining
arcs and each of these curves intersect $\mathcal{C'}$ in $\alpha$
points.

\begin{figure}
\begin{center}
\includegraphics*[totalheight=4cm]{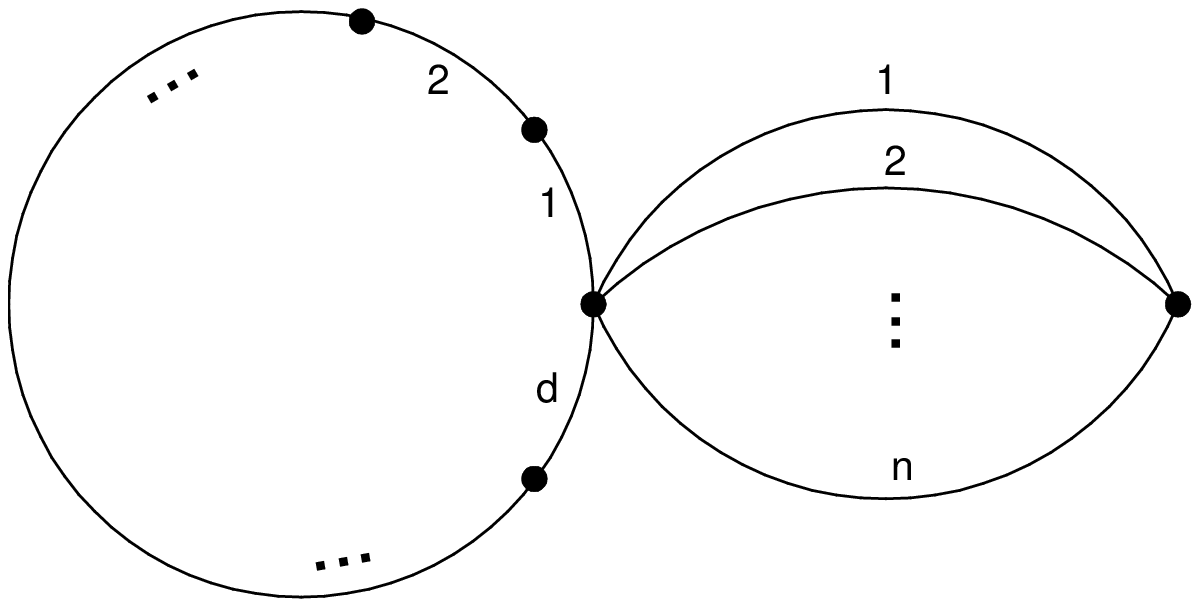}
\end{center}
\caption{} \label{due-ponti}
\end{figure}

The graph  $G(\mathcal{C}'')$ is the one depicted in Figure
\ref{due-ponti}, and each of its maximal tree $T$ consists of $d-1$
edges corresponding to curves of type A and one edge corresponding
to a curve of type B. So, the total number of vertices of
$\Gamma(H)$ that belong to curves corresponding to the edges of $T$
is $\alpha+(d-1)n/d$. By removing the curves corresponding to $T$
from $\mathcal{C''}$, we obtain a reduced Heegaard diagram which has
a  region, namely the upper one in Figure \ref{dun}, with at least $2n$
vertices. Indeed, except for sporadic cases, $2n$ is the maximal
number of vertices in a region. Anyway, the statement follows from
Theorem~\ref{prop}.
\end{proof}

An asymptotically equivalent estimate has been obtained
in~\cite{PV}, where a lower bound has been obtained in the
hyperbolic case (i.e. $\beta\ne 1, \alpha-1$) via volume arguments.
We give a lower bound for the remaining cases.

\begin{proposition}
Let $(n,s)\ne(3,1),(3,2)$ if $\alpha=2$. We have
$$c(C_{n,s}(\alpha,1)) = c(C_{n,s}(\alpha,\alpha-1)) \geqslant $$
$$
\geqslant 2 \log_5 \left(M\left(\frac{nm}{hD}\right)^m\left(
\frac{\alpha M}{2D}\right)^{M-1}\right) + D - M - m
$$
where $D=\gcd(n,\frac{\alpha}{2}(s-1))$, $M=\gcd(n,s-1)$,
$h=\gcd(n,s)$ and $m=\gcd(D,h)$.
\end{proposition}

\begin{proof}
Since $\mathbf b(\alpha,\alpha-1)$ is the mirror image of $\mathbf
b(\alpha,1)$ then $C_{n,s}(\alpha,\alpha-1)\cong C_{n,s}(\alpha,1)$.
Moreover $\mathbf b(\alpha,1)$ is the 2-component torus link of type
$(\alpha,2)$. So $c(C_{n,s}(\alpha,1))$ is a Seifert manifold and
then it is irreducible. The first homology group is computed in
\cite{Mu}.  So, the statement follows  applying  Theorem 2.6.2 of \cite{M2}.
\end{proof}

\subsection{A class of cyclic branched coverings of theta graphs}

\begin{figure}
\begin{center}
\includegraphics*[totalheight=1.7cm]{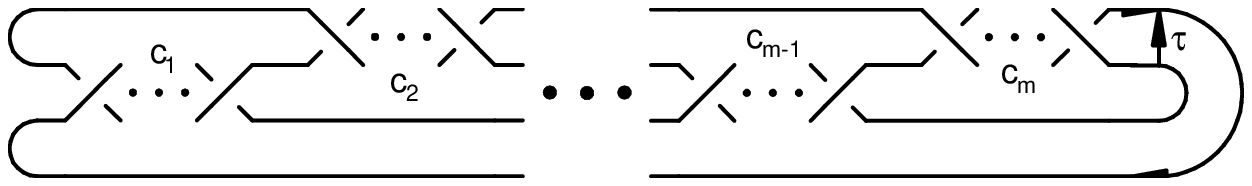}
\end{center}
\caption{The theta graph $\Theta(\alpha,\beta)$.} \label{tunnel}
\end{figure}

Let $\Theta(\alpha,\beta)$  be the theta graph in $\S^3$ obtained
from a two bridge knot of type $(\alpha,\beta)$ by adding a lower
tunnel $\tau$ as in Figure~\ref{tunnel}. Without loss of generality
we can assume that
$$
\frac{\alpha}{\beta} = c_1 + \frac{\displaystyle 1}{\displaystyle
c_2 + \, \cdots \, + \frac{\displaystyle 1}{\displaystyle c_{m-1} +
\frac{\displaystyle 1}{\displaystyle c_{m}}}} ,
$$
where $m>0$ and $c_1, \ldots, c_m$ can be taken as even integers
(see \cite[p.~26]{Kawa}).

Let $n>2$ and $s \in \mathbb Z_n - \{ 0, 1\}$, then we denote by $\Theta_{n,s}(\alpha, \beta)$ the
$n$-fold cyclic branched covering of $\Theta(\alpha,\beta)$ having
monodromy $\omega(m_1)=1$, $\omega(m_2)=s$ and $\omega(m_3)=s-1$,
 where $m_3$ is a meridian loop
around the tunnel and $m_1,m_2$ are meridian loops around the other
two edges of the graph, according to the orientations depicted in
Figure~\ref{tunnel}. By result of \cite{Mu2},
$\Theta_{n,s}(\alpha,\beta)$ is a pseudomanifold with two singular
points whose links are both homeomorphic to a closed surface of
genus $(1+n-\gcd(n,s)-\gcd(n,s-1))/2$.

\begin{proposition} Let $\widehat{\Theta}_{n,s}(\alpha,\beta)$ be
the compact manifold obtained by removing regular neighborhoods of
the two singular points of $\Theta_{n,s}(\alpha, \beta)$,  then
$$
c(\widehat{\Theta}_{n,s}(\alpha,\beta)) \leqslant n(\alpha-1).
$$
\end{proposition}

\begin{proof}
It follows from a result of \cite{Mu2} that
$\widehat{\Theta}_{n,s}(\alpha,\beta)$ is homeomorphic to the
generalized Dunwoody manifold
$M(\beta,\alpha-2\beta,1,n,2\beta-\alpha,s)$. Thus we can use
Theorem~\ref{prop} to calculate $\widetilde{c}(H)$ in order to
obtain an upper bound for $c(\widehat{\Theta}_{n,s}(\alpha,\beta))$.

\begin{figure}[h]
\begin{center}
\includegraphics*[totalheight=3cm]{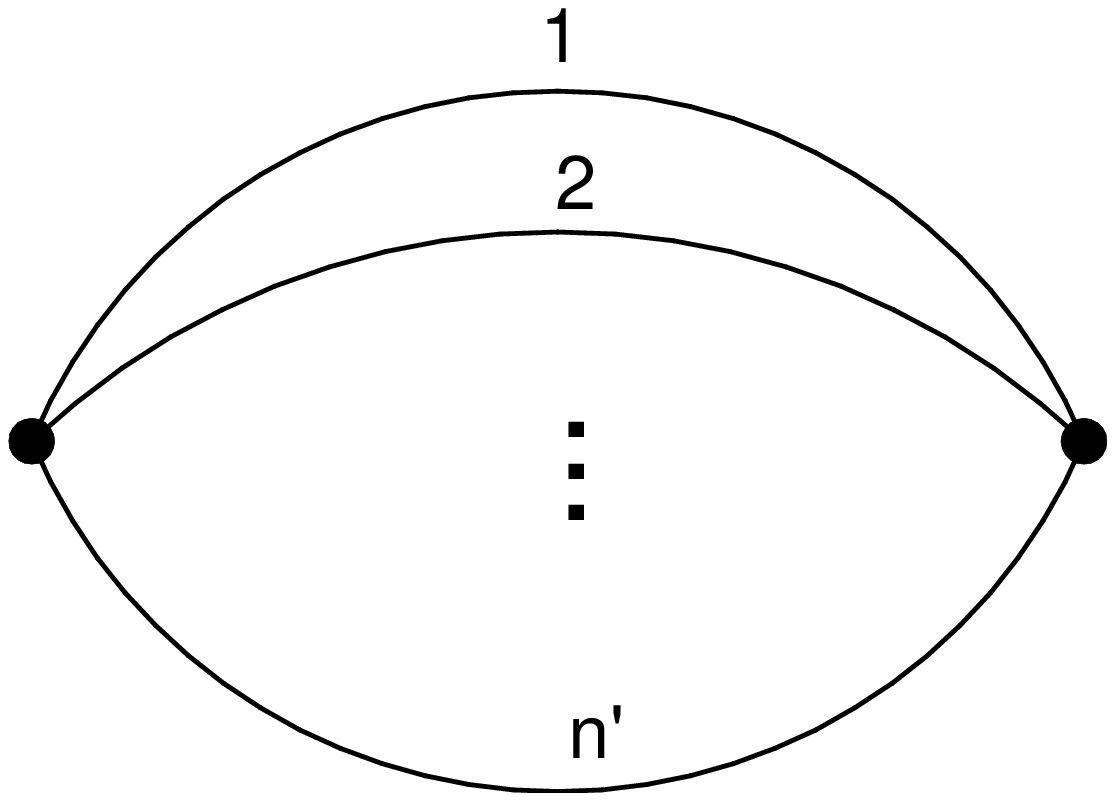}
\end{center}
\caption{} \label{theta}
\end{figure}

The system of curves $\mathcal C''$ of the Dunwoody diagram
$$
H=H(\beta,\alpha-2\beta,1,n,2\beta+1,s) =
(\Sigma_n,\mathcal{C'},\mathcal{C''})
$$
is reduced. Indeed, taking advantage of its symmetries, it is easy
to see that it  consists of $n'=\gcd(n,s)+\gcd(n,s-1)$ curves. More
precisely, $\gcd(n,s)$ curves arise from the ``radial'' arcs, while
the other $\gcd(n,s-1)$ curves  arise from the remaining arcs. The
graph $G(\mathcal{C}'')$ is the one depicted in Figure~\ref{theta},
where each vertex corresponds to a region of genus $(1+n-n')/2
> 0$. So $T$ consists of two isolated vertices and the system $\C''$
is already reduced. Since $\alpha$ is odd, then $\alpha \ne 0$ and
$\alpha - 2\beta \ne 0$. So, referring to Figure~\ref{dun}, the
region with the maximum number of vertices is always the upper one,
which has $2n$ vertices. The statement follows from
Theorem~\ref{prop}.
\end{proof}

\vspace{15 pt} {ALESSIA CATTABRIGA, Department of Mathematics,
University of Bologna, I-40126 Bologna, ITALY. E-mail:
cattabri@dm.unibo.it}

\vspace{15 pt} {MICHELE MULAZZANI, Department of Mathematics,
University of Bologna, I-40126 Bologna, ITALY. E-mail:
mulazza@dm.unibo.it}

\vspace{15 pt} {ANDREI VESNIN, Sobolev Institute of Mathematics,
Novosibirsk, \hbox{RUSSIA}. E-mail: vesnin@math.nsc.ru}

\end{document}